\documentclass[a4paper,11pt]{amsart}

\usepackage[T1]{fontenc}
\usepackage{mathrsfs, dsfont}
\usepackage{tikz}
\usetikzlibrary{calc}
\usetikzlibrary{intersections}
\usepackage{comment}
\usepackage{psfrag}

\definecolor{colore1}{RGB}{205,255,0}
\definecolor{colore2}{RGB}{255,255,0}
\definecolor{colore3}{RGB}{255,191,0}
\definecolor{coloreLinea1}{RGB}{0,0,255}
\definecolor{coloreLinea2}{RGB}{255,0,0}
\definecolor{coloreLinea3}{RGB}{148,0,105}

\tikzstyle{numeriLinea}=[black, pos=0.8, fill=white]

\usepackage{tikz}
\usepackage{comment}

\definecolor{colore1}{RGB}{205,255,0}
\definecolor{colore2}{RGB}{255,255,0}
\definecolor{colore3}{RGB}{255,191,0}
\definecolor{coloreLinea1}{RGB}{0,0,255}
\definecolor{coloreLinea2}{RGB}{255,0,0}
\definecolor{coloreLinea3}{RGB}{148,0,105}

\usepackage{amssymb,latexsym, amsthm, pifont}
\usepackage{hyperref}
\usepackage{cancel}
\usepackage[top=3.5cm,bottom=2.5cm,left=2.3cm,right=2.3cm,%
	bindingoffset=0cm]{geometry}
\usepackage{psfrag}
\usepackage{verbatim}

\usepackage{graphicx}
\usepackage{color}
\usepackage{nomencl} 

\theoremstyle{plain}
\newtheorem{lemma}{Lemma}[section]

\newtheorem{theorem}[lemma]{Theorem}

\theoremstyle{definition}

\theoremstyle{remark}

\numberwithin{equation}{section}

\makenomenclature
\RequirePackage{ifthen}
\renewcommand{\nomgroup}[1]{%
 \ifthenelse{\equal{#1}{A}}{\item[\bfseries{General mathematical symbols}]}{%
 \ifthenelse{\equal{#1}{B}}{\item[]\item[\bfseries{Common Symbols for Systems of Conservation Laws}]}{{%
 \ifthenelse{\equal{#1}{C}}{\item[]\item[\bfseries{Symbols introduced in the present paper}]}{}}}}}

\newcommand{\R}{\mathbb{R}}
\newcommand{\N}{\mathbb{N}}

\newcommand{\loc}{\text{\rm loc}}

\newcommand{\norm}[1]{\left\|#1\right\|}

\DeclareMathOperator{\TV}{\mathrm{TotVar}}
\DeclareMathOperator{\BV}{BV}
\DeclareMathOperator{\SBV}{SBV}

\newcounter{stepnb}
\newcounter{substepnb}

\title[Regularity results and counterexamples for conservation laws]{A note on regularity and failure of regularity for systems of conservation laws via Lagrangian formulation}
\author[L. Caravenna]{Laura Caravenna}
\address{L.C. Dipartimento di Matematica,
Universit\`a degli Studi di Padova,
Via Trieste 63, 35121 Padova, Italy}
\email{caravenna@math.unipd.it}

\begin{document}

\maketitle

{
\rightskip .85 cm
\leftskip .85 cm
\parindent 0 pt
\begin{footnotesize}

{\sc Abstract.}
The paper recalls two of the regularity results for Burgers' equation, and discusses what happens in the case of genuinely nonlinear, strictly hyperbolic systems of conservation laws. 
The first regularity result which is considered is Ole{\u\i}nik-Ambroso-De Lellis $\SBV$ estimate: it provides bounds on $\partial_{x} u$ when $u$ is an entropy solution of the Cauchy problem for Burgers' equation with $L^{\infty}$-data. Its extensions to the case of systems is then mentioned.
The second regularity result of debate is Schaeffer's theorem: entropy solutions to Burgers' equation with  $C^{k}$-data which are generic, in a Baire category sense, are piecewise smooth. The failure of the same regularity for general genuinely nonlinear systems is next described. The main focus of this paper is indeed including heuristically an original counterexample where a kind of stability of a shock pattern made by infinitely many shocks shows up, referring to~\cite{CaravennaSpinolo:counterex} for rigorous proofs.



\end{footnotesize}

}


\section{Introduction}
This paper makes a parallel among two regularity results which are very classical in the case of a single equation with uniformly convex flux
\begin{equation}
\label{E:Burgers}
\partial_{t}u +\partial_{x} f(u)=0,
\qquad f''(z)>0, \quad u:\R^{+}\times\R\to\R,
\end{equation}
but which have been discussed only in recent years for systems
\begin{equation}
\label{E:Burgers}
\partial_{t}U +\partial_{x} F(U)=0,
\qquad F\in C^{\infty}(\R^{N};\R^{N}), \quad U:\R^{+}\times\R\to\R^{N}.
\end{equation}
The analysis, via the Lagrangian formulation, concerns systems satisfying the assumptions of:
\begin{itemize}
\item Strict hyperbolicity: the Jacobian matrix $JF$ has real distinct eigenvalues $\lambda_{1},\dots,\lambda_{N}$.
\item Genuine nonlinearity: $ \nabla \lambda_i \cdot  \vec r_i  \ge c > 0 $ with $ r_i$ smooth $i$-th right eigenvectors of $JF$, for all $i$.
\end{itemize}

Strict hyperbolicity allows for a well-posedness theory, which was built with great contributions starting by Lax~\cite{Lax} and Glimm~\cite{Glimm}, for initial data having locally small total variation.
Good references are~\cite{Bre, Dafermos, HoldenRisebro}.
In particular, there is a lack of uniqueness of distributional solutions and the `admissible' solution is also selected as limit of wave-front tracking approximations.

Genuine nonlinearity is a uniform convexity assumption of the primitive of the restriction of $\lambda_{i}$ on the integral curves of $\vec r_i  $.
It generalizes the convexity assumption on $f$. 
Nonlinearity has regularizing effects, some holding also for multi-D equations~\cite{ChenFrid}: two of them are briefly recalled below. 

\vskip.5\baselineskip
Let us first briefly highlight few achievements in the case of a single convex equation.
Starting with a datum $u_{0}\in L^{\infty}$, Ole{\u\i}nik~\cite{Oleinik} shows that for every $a<b$
\[
u(t,b)-u(t,a)\leq \frac{b-a}{t c}, \qquad c=\inf f''.
\]
This upper estimate only allows for decreasing jumps and it implies that $L^{\infty}$ data are instantaneously regularized to functions of locally bounded variation ($\BV_{\loc}$).
The same lower estimate cannot hold, otherwise one would get the Lipschitz continuity at time $t$.
Ambrosio-De Lellis~\cite{AmbrosioDeLellis} improved the regularity to the space of \emph{special} functions of bounded variation ($\SBV_{\loc}$) via a lower estimate: denoting by $X(t,x)$ the map moving points forward along characteristic curves, it is of the kind
\[
u(t,b)-u(t,a)\geq -2\frac{d-c}{tc} ,
\]
where $d=\inf\{x\ : X(t,x)> b\}$ and $c=\sup\{x\ : X(t,x)< a\}$. See~\cite{AFP} for $\BV$ and $\SBV$ spaces.

\vskip.5\baselineskip
There are other interesting decay estimates which are not reviewed here, for passing to a different regularization effect by Schaeffer~\cite{Schaeffer}.
Schaeffer's regularity theorem, considering also its improvement by Dafermos~\cite{Daf}, sates the generic piecewise smooth regularity of the entropy solutions of the Cauchy problem of~\eqref{E:Burgers} for initial data in the Schartz space $\mathfrak S$ of $C^{\infty}$ and rapidly decreasing functions.

\begin{theorem}[Schaeffer, 1973--Dafermos, 1985]
There is an open, dense subset of $\mathfrak S$ such that for initial data belonging to this subset the admissible solution to the Cauchy problem of~\eqref{E:Burgers} is piecewise smooth: there are finitely many smooth shocks parametrized by the time, $u$ is smooth in the complementary of the shocks, it has a continuous limit from both sides of any shock and this limit is smooth at interior points of the curves of shocks.
\end{theorem}
According to Schaeffer's opinion, ``the proof of this theorem is perhaps more interesting than the result itself, which is not too surprising''. His approach is strictly related to the uniform convexity of the flux function, but with the method of generalized characteristics several later improvements have been reached for a single equation, just mention~\cite{BiaYuScalar} within BV.
Let us turn instead the attention to the case of systems.

 \vskip.5\baselineskip
Decay estimates can be generalized to genuinely nonlinear systems, see e.g.~\cite{BressanColombo,BressanGoatin,BressanYang, BiaCar, GlimmLax,Liu_decay, ChristoforouTrivisa}.
As for the case of a single equation, they yield that initial data small in $\BV$ are regularized to the space of functions of special bounded variation~\cite{AnconaNguyen,BiaCar,BianchiniYu}.
See moreover~\cite{Dafermos_sbv} beyond genuinely nonlinearity, but for Riemann problems.
The analyses of systems have a very different and much more complicated proof, but it is not just technique as they present new behaviors. However, the analysis is still based on the Lagragrangian approach, at an approximated level. 
One would likely have easier proofs if a Lagrangian description was already available directly on the admissible solution, instead of having it only at an approximated level.

What about Schaeffer's regularity theorem in the case of systems? The main point of this paper is showing that Schaeffer's theorem cannot be generalized to the class of genuinely nonlinear systems.
The original partial counterxample that was presented at HYP is here outlined, as the one in~\cite{CaravennaSpinolo:counterex}.
Even if shocks for conservation laws are stable objects, Schaeffer \emph{proves} that patterns with infinitely many shocks are not stable, for a single convex equation.
In the case of genuinely nonlinear systems, instead, the examples below show that shock patterns with infinitely many shocks similar to the ones in~\cite{BaJ} may be robust. 

{Why should one bother of extending Schaeffer's regularity to the case of systems?} Positive extensions are relevant for their consequences in analysis and in numerical analysis: just compare the different complexity of the proof of convergence of vanishing viscosity approximations by Goodman \and Xin~\cite{GoodmanXin}, for piecewise smooth solutions, with the general one by Bianchini \and Bressan~\cite{BiaBre}. Negative answers widen as well our knowledge of the foundations of models that are close to physics, because they show the possibility that infinite shock patterns are stable.

Observe finally that the mechanism creating the `stable' shock pattern is built up locally: \emph{the total variation being small in this case is not a limit, but a plus}.
The total variation is small \emph{locally} in the region that is considered and where the shock pattern with infinitely many shocks appears, but it may be bigger globally. I shall also emphasize that a technical tool at the base of our proof, namely the fine convergence properties of the wave-front tracking approximations, was not present two decades ago.

\section{The case of Genuinely Nonlinear Strictly Hyperbolic systems}

For the case of systems it is not a priori evident if the generic piecewise smooth regularity holds or not. In the case of systems, on the one hand the number of waves may grow exponentially, which would suggest a negative answer; on the other hand, they are not necessarily shock waves, which carry discontinuities: the growing number of rarefaction waves, which are continuous, could enforce cancellations. 

The answer is not. Below there is a description of two counterexamples where the shock structure with infinitely many shocks posses kinds of stability properties. The counterexamples rule out the possibility that generic piecewise smooth regularity holds for the class of genuinely nonlinear systems. 
Our aim here is to provide precise statements of both the different examples and an accessible heuristic description of the first counterexample, presented at HYP2014.
There is also a map for the full proof of this first original counterexample which refers to precise definitions and computations published in~\cite{CaravennaSpinolo:counterex} relatively to the second one: the reader who is looking for a full understanding of the detail of the proof shall read first that paper.

\subsection{Prelimiaries}
\label{Ss:prelimiaries}
For both the examples, consider the system provided by Baiti and Jenssen~\cite{BaJ}:
\begin{equation}
\label{E:BaitiJenssen}
\partial_{t}\left(\begin{matrix} u\\ v\\ w
	\end{matrix}\right) 
+
\partial_{x}\left(\begin{matrix}
	\displaystyle{ 4 \big[ (v-1) u - w \big] + 2 \eta \left[ u w-  u^{2} (v-1 )\right] \phantom{\int}} \\
	v^2 \\
	\displaystyle{4 \Big\{ {v (v-2)u} - (v-1) w \Big\}
	 + \eta \left[w^2-u^2 (v-2) v\right] } \\
	\end{matrix}\right) 
=0,
\end{equation}
where $u$, $v$, $w$ are real functions defined in $\R^{+}\times\R$ and $\eta \in\left]0,\frac{1}{4}\right[$. This system is strictly hyperbolic:
\begin{align*}
&\lambda_{1}(U)=2 \eta \big[w-(v-2)u \big]-4 && \in [-6, -5/2] \\
& \lambda_{2}(U)=2v &&\in [-2,2] \\
& \lambda_{3}(U)=2 \eta \big[ w-v u \big]+4 &&\in [3,5]
\end{align*}
when $|U|<1$ and $0<\eta<1/4$.
It is also genuinely nonlinear because
\[
\nabla\lambda_{1}(U)\cdot \vec r_{1}(U)=4 \eta =-\nabla\lambda_{3}(U)\cdot \vec r_{3}(U),
\qquad
\nabla\lambda_{2}(U)\cdot \vec r_{2}(U) = 2,
\]
where above $U=(u,v,w)^{T}$ and there is a fixed set of right eigenvectors of $JF$:
\begin{equation*}
\vec r_{1}(U)=\begin{pmatrix}1 \\ 0 \\ v\end{pmatrix},
\quad 
\vec r_{3}(U)=\begin{pmatrix}1 \\ 0 \\ v-2 \end{pmatrix},
\quad\text{$\vec r_{2}$ has second component $1$.}
\end{equation*}

Given an admissible distributional solution $U$ of~\eqref{E:BaitiJenssen}, a \emph{k-shock} is a jump discontinuity of $U$ which propagates in space-time with speed in the range of $\lambda_{k}$; this speed is called \emph{speed} of the $k$-shock, and $k=1,2,3$ is called \emph{family} of the wave.
Concerning the special system~\eqref{E:BaitiJenssen}, the \emph{strength} of an $i$-shock from the left state $U^{-}$ to the right state $U^{+}$ is
\begin{itemize}
\item for $i=1,3$ the modulus of the constant $s$ of proportionality $U^{+}-U^{-}=s \vec r_{i}(U^{-})$,
\item for $i=2$ the modulus of the difference of the second components $|s|=|v^{+}-v^{-}|$.
\end{itemize}
\emph{Rarefaction waves} are continuous waves. For Burgers' equation, corresponding to $f(z)=z^{2}/2$, for example they arise at Riemann problems when the left state $u^{-}$ is lower than the right state $u^{+}$: the admissible solution $u(t,x)=u^{-}\chi_{\{x\leq tu^{-}\}}+\frac{x}{t} \chi_{\{tu^{-}< x< tu^{+}\}}+u^{+}\chi_{\{x\geq tu^{+}\}}$ is Lipschitz continuous with a wave fan emanating from the origin. 
In the wave-front tracking approximation, which is a piecewise constant approximation converging to the admissible solution, rarefactions are discretized into small jumps. Let us adopt notations above also for strengths and speed of this discretization, see more on~\cite{Bre}.
In this exposition, there is the effort to hide the technicality related to the wave-front tracking approximations.
Finally, let us recall that if $u$ is the above rarefaction wave for Burgers, then $v(t,x)=u(1-t,-x)$ is a Lipschitz continuous solution for $0\leq t<1$ which develops at time $1$ the jump $u^{+}$ (left) $u^{-}$ (right): term $v$ \emph{compression wave}. It can be also generalized to systems: $k$-compression waves are Lipschitz solutions of~\eqref{E:BaitiJenssen} producing $k$-shocks.

A good feature of system~\eqref{E:BaitiJenssen} is the following: there are reasonable sufficient conditions on the left and right states of a Riemann problem which ensure that the admissible solution is obtained patching together only, in the order, 1-, 2-, and 3-shocks, and no rarefaction.
Moreover, lower bounds on the strengths of the shocks can be provided.
They were first deduced in~\cite{BaJ} for establishing that all the Riemann problems that they encounter, both at the initial time and later when resolving interactions of waves, were producing shocks sufficiently strong, and only shocks, as described below.
Their analysis is extended in~\cite{CaravennaSpinolo} with new estimates and new cases needed for the stability properties of the counterexample~\cite{CaravennaSpinolo:counterex}.
Another simplifying feature of system~\eqref{E:BaitiJenssen} is that when 1- and 3-waves interact among themselves the nonlinear interaction only generates the same kind of 1- and 3-waves with equal strength, and possibly different speed.

\subsection{An infinite shock pattern}
\label{Ss:infiniteshockpattern}
Before discussing stability properties, let us see at all a shock pattern consisting of countably may shocks: let us briefly term it \emph{infinite shock pattern}. This is along the lines of~\cite{BaJ}, but in the context of small total variation and with precise values different from there. The presentation here is a bit more complex than necessary because it will be a reference for \S~\ref{Ss:counterexample}.
Good values of the parameters both for this section and for the analysis in \S~\ref{Ss:counterexample} are~\cite{CaravennaSpinolo}:
\begin{align*}
 &q= 20, & \eta : = \varepsilon^2, &&  \omega := \varepsilon^3, && r := \varepsilon^{10}/4, 
 &&\widetilde T={40}/{\varepsilon^{3}} ,
 &&\rho := 12 \widetilde T + 40 ,  &&a=q+7
\end{align*}
and $\varepsilon>0$ is a sufficiently small parameter.

Consider the Cauchy problem for~\eqref{E:BaitiJenssen}.
For small $s>0$ and $|U^{-}|<1$, denote\footnote{This notation for the shock waves slightly differs from the one in~\cite{CaravennaSpinolo:counterex}, but it is enough here.} by $S_{i}[s,U^{-}]$ the right value of an admissible $k$-shock having left value $U^{-}$ and strength $s$.
The \emph{piecewise constant} initial datum
\begin{gather*}
Z(x)=U_{I}\chi_{\{x<-q\}} + U_{I\!I}\chi_{\{-q<x<q\}} +U_{I\!I\!I}\chi_{\{x>q\}} ,\\
 U_{I\!I}=S_{3}[\omega,S_{2}[\omega,[S_{1}[\omega,U_{I}]]], \quad  U_{I\!I\!I}=S_{3}[\omega,S_{2}[\omega,[S_{1}[\omega,U_{I\!I}]]], \quad U_{I}=(\varepsilon,\omega,-\varepsilon)
\end{gather*}
generates, for small times, from left to right, the following waves:
\begin{enumerate}
\item\label{item:1lex1} a 1-shock close to $x=-q$ of order $\omega'$ between $\omega\sqrt{\varepsilon}$ and $\omega$ with speed in $[-6, -5/2]$,
\item\label{item:2lex1} a 2-shock $J_{\ell}$ of the order $\omega$ close to $x=-q$ with {speed} close to $\omega$,
\item\label{item:3lex1} a 3-shock $R_{0}$ close to $x=-q$ of order between $\omega\sqrt{\varepsilon}$ and $\omega$ with speed in $[ 3,5]$ and 
\item\label{item:1rex1} a 1-shock $S_{0}$ close to $x=q$ of order between $\omega\sqrt{\varepsilon}$ and $\omega$ with speed in $[-6, -5/2]$,
\item\label{item:2rex1} a 2-shock $J_{r}$ of the order $\omega$ close to $x=q$ with speed close to $-\omega$,
\item\label{item:3rex1} a 3-shock close to $x=q$ of order $\omega''$ between $\omega\sqrt{\varepsilon}$ and $\omega$ with speed in $[ 3,5]$, 
\end{enumerate}
See notations in Figure~\ref{fig:esempio} a bit after the initial time, even if the picture includes also other elements that are introduced below.
Claim: a {piecewise constant} initial datum whose short time evolution is given by \ref{item:1lex1})-\ref{item:3rex1}) develops a shock pattern with infinitely many shocks~\cite{BaJ}. The rough picture is the following:
\begin{itemize}
\item The 2-shocks $J_{\ell}$ and $J_{r}$ interact at some time $t=\widetilde T$ close to $t=q/\omega$.
\item The 1- and 3-shocks $S_{0}$ and $R_{0}$ interact. They then prosecute with the same strength, but with different speed due to the nonlinearity of the interaction.
\item The 1-shock $S_{0}$ interacts with the left 2-shock $J_{\ell}$. After this, the left 2-shock $J_{\ell}$ proceeds with the same speed and strength.
The transmitted 1-shock leaves the triangle with edges drawn by $J_{\ell}$ and $J_{r}$. A new 3-shock $R_{1}$ is generated, with strength of the order $\omega\omega'$.
\item The 3-shock $R_{0}$ interacts with the right 3-shock $J_{r}$. After this, the left 2-shock $J_{r}$ proceeds with the same speed and strength.
The transmitted 3-shock leaves the triangle with edges drawn by $J_{\ell}$ and $J_{r}$. A new 1-shock $S_{1}$ is generated, with strength of the order $\omega\omega''$.
\end{itemize}
After $j$ iterations, one can still repeat the last points because, owing to geometric constraint given by the restriction on speeds of waves of different families, it is before the time $\widetilde T$ of 2-2 interaction:
\begin{itemize}
\item The 1- and 3-shocks $S_{j}$ and $R_{j}$ interact, and they proceed up to the interaction with $J_{r} $ and $J_{\ell}$, respectively. The transmitted 1- and 3-shocks leave the triangle with edges drawn by $J_{\ell}$ and $J_{r}$. The reflected 3- and 1-shocks $R_{j+1}$ and $S_{j+1}$ are of orders $\omega^{j+1}\omega'$ and $\omega^{j+1}\omega''$. 
\end{itemize}
In Baiti and Jenssen construction the strengths of $S_{j}$ and $R_{j}$ were increasingly amplified instead of dumped, mechanism which blows-up the $L^{\infty}$-norm.

In the counterexamples below we rather construct two initial data which are
\begin{enumerate}
\item Smooth. With the same generality, they can be chosen compactly supported.
\item \label{item:smallnessdata} Within the framework of the well-posedness theory for systems, see for instance~\cite{Bre}. 
\item Robust in developing countably many shocks locally in space-time, if one adds perturbations that are small in suitable spaces. The two examples differ for the space of perturbations.
\end{enumerate}
Observe how \ref{item:smallnessdata}) implies that if there are countably many shocks at a same time in a compact set they must be increasingly smaller, because their strength must sum up to a value which is uniformly bounded: these shocks are therefore in principle fragile in surviving perturbations, as pointed out next.

\subsection{First counterexample}
\label{Ss:counterexample}
Hyperbolic equations with bounded initial data have finite speed of propagation. The discussion can~\cite[Remark~5.1]{CaravennaSpinolo} thus be confined to a rectangle $[-\rho,\rho]\times[0,2\widetilde T]$ in space-time where the infinite shock pattern shows up, and which is not affected by the values of the initial datum out of $[-\rho,\rho]$.

\begin{figure}
\begin{center}
\begin{picture}(0,0)%
\includegraphics{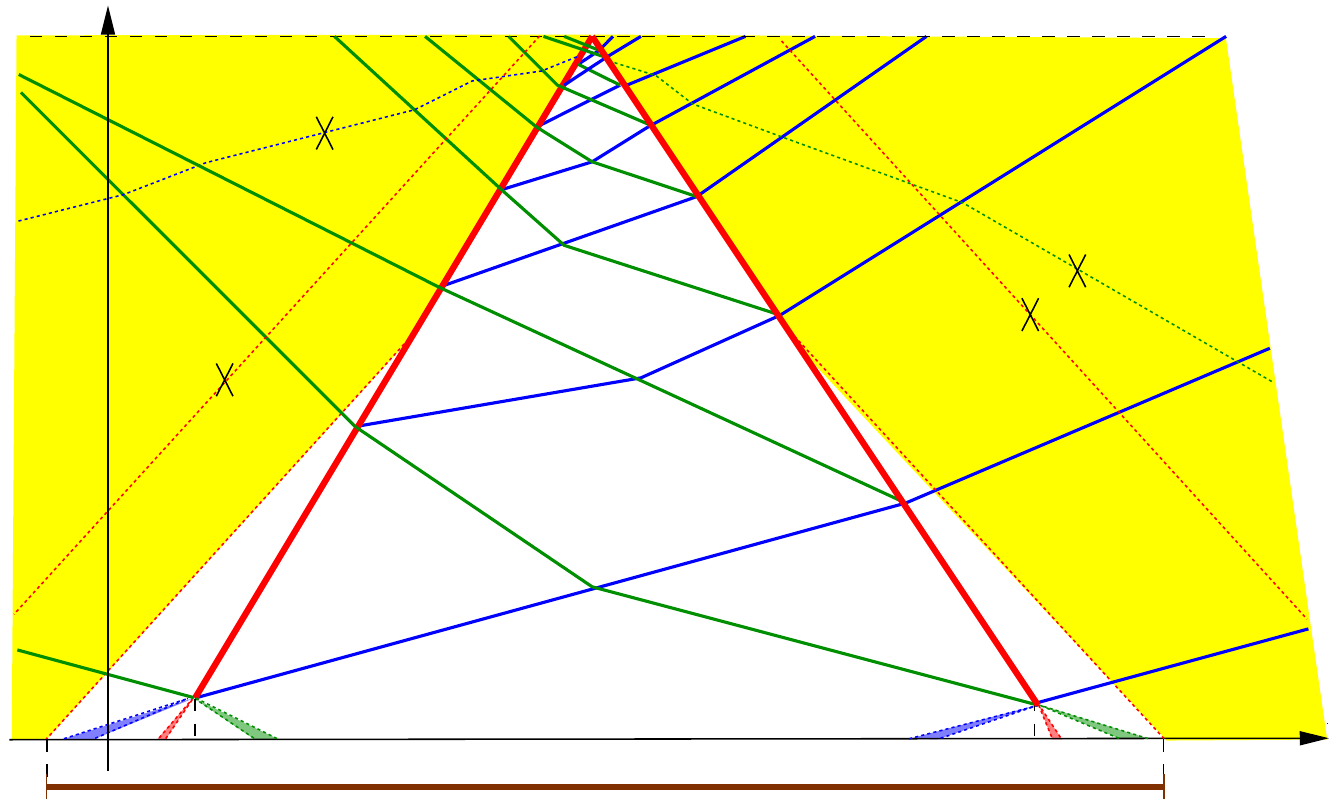}%
\end{picture}%
\setlength{\unitlength}{3947sp}%
\begingroup\makeatletter\ifx\SetFigFont\undefined%
\gdef\SetFigFont#1#2#3#4#5{%
  \reset@font\fontsize{#1}{#2pt}%
  \fontfamily{#3}\fontseries{#4}\fontshape{#5}%
  \selectfont}%
\fi\endgroup%
\begin{picture}(6391,3868)(83,-8867)
\put(1798,-8284){\makebox(0,0)[lb]{\smash{{{}{\color[rgb]{0,0,1}$\omega''$}%
}}}}
\put(1798,-8000){\makebox(0,0)[lb]{\smash{{{}{\color[rgb]{0,0,1}$R_{0}$}%
}}}}
\put(3300,-7609){\makebox(0,0)[lb]{\smash{{{}{\color[rgb]{0,0,1}$\omega''$}%
}}}}
\put(2678,-7609){\makebox(0,0)[lb]{\smash{{{}{\color[rgb]{0,.5,0}$\omega'$}%
}}}}
\put(3511,-6953){\makebox(0,0)[lb]{\smash{{{}{\color[rgb]{0,.5,0}$\omega \omega''$}%
}}}}
\put(3511,-7200){\makebox(0,0)[lb]{\smash{{{}{\color[rgb]{0,.5,0}$S_{1}$}%
}}}}
\put(3900,-8278){\makebox(0,0)[lb]{\smash{{{}{\color[rgb]{0,.5,0}$\omega'$}%
}}}}
\put(3900,-8000){\makebox(0,0)[lb]{\smash{{{}{\color[rgb]{0,.5,0}$S_{0}$}%
}}}}
\put(2202,-6857){\makebox(0,0)[lb]{\smash{{{}{\color[rgb]{0,0,1}$\omega \omega'$}%
}}}}
\put(2202,-7120){\makebox(0,0)[lb]{\smash{{{}{\color[rgb]{0,0,1}$R_{1}$}%
}}}}
\put(4400,-7800){\makebox(0,0)[lb]{\smash{{{}{\color[rgb]{1,0,0}$J_{r}$}%
}}}}
\put(4743,-7800){\makebox(0,0)[lb]{\smash{{{}{\color[rgb]{1,0,0}$\omega$}%
}}}}
\put(1437,-7841){\makebox(0,0)[lb]{\smash{{{}{\color[rgb]{1,0,0}$J_{\ell}$}%
}}}}
\put(1137,-7841){\makebox(0,0)[lb]{\smash{{{}{\color[rgb]{1,0,0}$\omega$}%
}}}}
\put(2438,-6400){\makebox(0,0)[lb]{\smash{{{}{\color[rgb]{0,0,1}$\omega^2 \omega''$}%
}}}}
\put(3113,-6228){\makebox(0,0)[lb]{\smash{{{}{\color[rgb]{0,.5, 0}$\omega^2 \omega'$}%
}}}}
\put(6100,-8461){\makebox(0,0)[lb]{\smash{{{}{\color[rgb]{0,0,0}$x$}%
}}}}
\put(646,-5100){\makebox(0,0)[lb]{\smash{{{}{\color[rgb]{0,0,0}$t$}%
}}}}
\put(346,-5211){\makebox(0,0)[lb]{\smash{{{}{\color[rgb]{0,0,0}$\widetilde T$}%
}}}}
\put(5750,-8770){\makebox(0,0)[lb]{\smash{{{}{\color[rgb]{0,0,0}$b$}%
}}}}
\put( 160,-8770){\makebox(0,0)[lb]{\smash{{{}{\color[rgb]{0,0,0}$a$}%
}}}}
\put( 900,-8650){\makebox(0,0)[lb]{\smash{{{}{\color[rgb]{0,0,0}$-q$}%
}}}}
\put( 5000,-8650){\makebox(0,0)[lb]{\smash{{{}{\color[rgb]{0,0,0}$q$}%
}}}}
\end{picture}%
\caption{The developing of an infinite shock pattern, see Theorem~\ref{T:CSa}. Compression waves in $]-a,a[$ generate two 2-shocks of order $\omega$, and two couples of 1- and 3-shocks of order $\omega'$. Preventing waves at initial time outside the interval $]-a,a[$, then 1- and 2-waves cannot be present in the left yellow region, nor 2,3-waves in the right one. Close to $- a$ and $a$ there are 2-characteristics which merge with $J_{\ell}$ and $J_{r}$ by the time $t=8/\omega$.
}
\label{fig:esempio}
\end{center}
\end{figure}
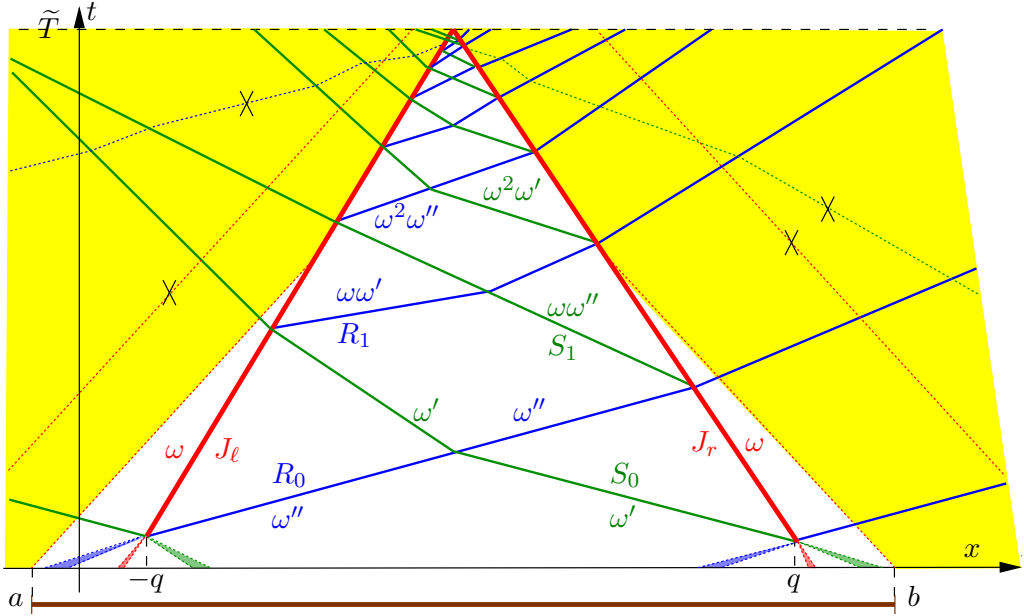

Let us state the first result, commenting it later. See Figure~\ref{fig:esempio}, that will be better described below highlighting the construction.
The counterexample in~\cite{CaravennaSpinolo:counterex}, see \S~\ref{Ss:secondschaeffer}, is stronger in the sense that here we only allow perturbations confined in a given interval. As well, the stability of the shock structure of this first example is qualitatively interesting because there is a different behavior than there, allowing rarefaction waves.
Moreover, Theorem~\ref{T:CSa} already prevents straight extensions of Schaeffer's regularity: for a single convex equation this statement would be false. Indeed, Schaeffer's regularity theorem also guarantees that the infinite shock pattern is not stable if one perturbs the initial datum in an interval containing the support of the derivatives.

\begin{theorem}[C., Spinolo-2014]
\label{T:CSa}
Consider the Cauchy problem of~\eqref{E:BaitiJenssen}.
There exists a $C^{\infty}$ initial datum $\widetilde U$ with small total variation in $]-a,a[$, constant outside $]-a,a[$, and a constant $r>0$ such that every initial data belonging to the nonempty ball 
\[
\mathfrak B = \widetilde U+ \left\{ Z \text{ compactly supported in $]-a,a[$ t.c.~} \TV(Z)<r \right\}
\]
develops infinitely many shocks in a compact set: it contains a shock pattern like in Figure~\ref{fig:esempio}.
\end{theorem}

The proof of Theorem~\ref{T:CSa} is now pictured with the main insights, like the definition of $\widetilde U$ and why cutting down perturbations out of $]-a,a[$.
The reader who is interested in the detail can reconstruct it from~\cite{CaravennaSpinolo:counterex}.

\subsubsection*{{\sc Step 0}: Preliminary considerations} In~\cite{BaJ} the solution was piecewise smooth and therefore the analysis could be performed exactly. With a smooth initial data direct computations are not that clear, presently. Theorems~\ref{T:CSa} and~\ref{T:CSb} rely instead on the wave-front tracking approximation and on its fine properties of convergence \cite{Bre,BressanLeFloch}.
What one cares to prove is that the wave-front tracking approximations of the admissible solution $U$ exhibit the infinite shock pattern with estimates from below on the {strengths} of the shocks uniformly in the approximation parameter.

For Theorem~\ref{T:CSa}, one therefore needs to understand here why the infinite shock pattern is present starting with a piecewise constant function $Z_{0}$, with $\TV(Z_{0}-\widetilde U)<r$, which evolves via wave-front tracking.
Having suitable lower estimates on the jumps, \cite[\S~5.7]{CaravennaSpinolo:counterex} presents a similar and detailed discussion on the technical limiting procedure.

\subsubsection*{{\sc Step 1/2:} Construction of $\widetilde U$}
The rough idea of smooth initial data developing infinitely many shocks is the following: providing a smooth, small initial datum whose profile at time $t=1$ resembles the profile that Baiti \and Jenssen had at initial time, which then develops a shock pattern like in \S~\ref{Ss:infiniteshockpattern}. The analysis is later different, dealing with small shocks rather than big shocks, having perturbations and being in the context of global existence; nevertheless, one can recognize among the many waves a pattern of infinitely many shocks. 

More precisely, produce first a Lipschitz continuous initial profile $V$ patching at time $t=0$ two blocks each made of 3-, 2-, 1-compression waves that collapse at $t=1$ into the 1-, 2-, 3-shocks described by \ref{item:1lex1})-\ref{item:3rex1}) of \S~\ref{Ss:infiniteshockpattern}. See (4.17) in~\cite{CaravennaSpinolo:counterex} for the explicit expression. A suitable mollification of $V$ eventually defines a desired $C^{\infty}$-initial profile $\widetilde U$, satisfying $\TV(\widetilde U-V)<r$.
Note how \S~\ref{Ss:infiniteshockpattern} tells us that $V$ develops the infinite shock pattern, and that the next steps deal with initial data $Z_{0}$ which are perturbations of $V$.


\subsubsection*{{\sc Step 1}: Initial formation of the `big' shocks.} One proves that at some time $\bar t$ close to $t=1$ the profile $Z(\bar t)$ evolved from $Z_{0}$ via wave-front tracking is close to the profile $V(t=1)$: this means that among the many waves there are six shocks satisfying \ref{item:1lex1})-\ref{item:3rex1}) of \S~\ref{Ss:infiniteshockpattern} together with
\begin{enumerate}\addtocounter{enumi}{6}
\item\label{item:o5ex1} possibly other waves, where the total strength of 1- and 3-shocks is at most of order $\omega$, while the total strength of the remaining waves is at most of the order $r$.
\end{enumerate}
The analysis rely on the subdivision of waves present at initial times into different subgroups:
\begin{itemize}
\item[A:] i-shocks arising at time $t=0$ in the region associated to the i-compression waves;
\item[B1-B2:] other shock waves arising at time $t=0$ and rarefaction waves arising at time $t=0$;
\item[C1-C2:] shock waves arising at later times and rarefaction waves arising at later times.
\end{itemize}
The evolution of the waves of each group is then tracked and their strengths and speed are carefully estimated by interaction estimates: the `big shocks' that we are interested in arise from Group A, while the rest is suitably controlled.
Neglecting discussions on the type of perturbation waves \ref{item:o5ex1}), which here may also be rarefactions, estimates leading to the above points are difficult but extremely close to the analysis in~\cite{CaravennaSpinolo:counterex}: Lemma~5.4 for~\ref{item:2lex1})-\ref{item:2rex1}), Lemma~5.3 for~\ref{item:1lex1}), \ref{item:3lex1}), \ref{item:1rex1}), \ref{item:3rex1}) and Lemma 5.2 for~\ref{item:o5ex1}).
To be fair, for achieving more straightly the above estimates, compression waves in Figure~\ref{fig:esempio} are slightly separated.
Estimates in~\ref{item:1lex1}), \ref{item:3lex1}), \ref{item:1rex1}), \ref{item:3rex1}) are not sharp, but they suffice for our aim.


\subsubsection*{{\sc Step 2}: Formation of the infinite shock pattern in the wave-front tracking approximations.} The reader finds here a sketch of the developing of the infinite shock pattern starting from an initial profile whose short time evolution consists of the waves~\ref{item:1lex1})-\ref{item:o5ex1}) above, together with estimates from below on the strengths of relevant shocks which are uniform in the wave-front tracking parameter, parameter in this exposition phantom.
Since {\sc Step 1} tells us that at time $\bar t$ there is such profile, this achieves the goal stated in {\sc Step 0}.

Recall from \S~\ref{Ss:infiniteshockpattern} that the infinite shock pattern is formed in absence of the perturbations waves of \ref{item:o5ex1}).
Interactions however are nonlinear: one would not otherwise have the blow-up in $L^{\infty}$ when initial data are large. The presence of perturbation waves makes therefore the analysis subtle, and different, even for the piecewise-constant wave-front tracking approximations.

\vskip.5\baselineskip
When only shock waves are present, the unperturbed analysis of \S~\ref{Ss:infiniteshockpattern} survives: what underpins the formation of the infinite shock structure in this situation is that the reflected waves $R_{j}$, $S_{j}$
\begin{itemize}
\item increase in strength if interacting with other waves of the same family, which are shocks;
\item interacting with waves of different families might be dumped by factors whose product is close to $1$, owing to interaction estimates and to the estimates on the lower order of perturbations at all times~\cite[Lemma 5.2]{CaravennaSpinolo:counterex};
\item at interactions with the 2-shocks $J_{\ell}$, $J_{r}$ still generate shocks satisfying lower estimates.
\end{itemize}
The full proof is in \cite[\S~5.7]{CaravennaSpinolo:counterex}. What happens instead if rarefaction waves are present?

\subsubsection*{Why perturbations of $V$, in this example, must be confined in the interval $]-a,a[$?}
Perturbations can easily cancel the infinite shock pattern if one of the following happens:
\begin{itemize}
\item There are 3-waves at times close to $t=\widetilde T$ at the left of $J_{\ell}$.
\item There are 1-waves at times close to $t=\widetilde T$ at the right of $J_{r}$.
\item There are 2-waves at times close to $t=\widetilde T$ at the left of $J_{\ell}$, or at the right of $J_{r}$, or both.
\end{itemize}
Indeed, being in the context of small total variation, the strengths of the shocks being formed constitue the terms of a convergent sum.
Assume there is a small perturbation of the initial datum that creates a 3-rarefaction which reaches the shock pattern close to the apex, but just before the time $\widetilde T$ of interaction among $J_{\ell}$ and $J_{r}$, see the curve erased with `X' in Figure~\ref{fig:esempio}.
Regardless of how small the 3-rarefaction is when it hits the left 2-shock $J_{\ell}$, it is possible that it reaches the shock pattern close enough to the vertex to be bigger than the 3-shocks $\{R_{j}\}_{k\in\N}$ which are being reflected and increasingly dumped: it is then possible that it erases one of them, and therefore the later ones, canceling the tail of the convergent sum of their strengths and obstructing the formation of the successive infinitely many shocks.
Other similar scenarios are possible if allowing one of the three possibilities above.
In this case, only \emph{finitely many} shocks would survive.

If perturbations of the initial data are small enough and confined in $]-a,a[$ one can show that, due to geometric constraints, this kind of cancellation is prevented: after time $t=8/\omega<\widetilde T/2$
\begin{itemize}
\item no 2- or 3-waves can be present on the left side of $J_{\ell}$ and
\item no 2- or 1-waves can be present on the right side of $J_{r}$. 
\end{itemize}
The first point is now naively explained as an example.
Suppose there is no wave initially present in $(-\infty,-a]$.
The 3-waves arising in $]-a,a[$ are faster than 2-waves due to the constraint on their speed recalled in \S~\ref{Ss:prelimiaries}: at time $t=8/\omega$ they are thus necessarily on the right side of $J_{\ell}$, after having interacted with other waves.
In this system, new 3-waves can be created only at interactions involving a 2-wave, and they stay then at the right of the 2-wave that generated it because of the speed constraint. If no 2-wave is present on the left side of $J_{\ell}$ after time $t=8/\omega$, then no 3-wave can be present after time $t=8/\omega$.
On the other hand, in this system no 2-wave is generated at positive times: due to speed constraints, one can show that the 2-waves possibly present at the left of $J_{\ell}$ at time $t=\bar t$ merge with $J_{\ell}$ by time $t=8/\omega$ because the speed of $J_{\ell}$ is about $\omega$ while the one of the 2-characteristic starting from $a$ is about $2\omega$.
This yields that no 2-wave is present after time $t=8/\omega$ on the left side of $J_{\ell}$, and consequently no 3-wave.

\subsubsection*{Why rarefactions which may arise in $]-a,a[$, or at later times, do not cancel the shock pattern?}
The idea of the main argument is the following.
Consider at time $t$ between $\bar t$ and $\widetilde T$
\begin{itemize}
\item the sum $\mathcal V_{13}(t)$ of the strengths of the 1- and 3-shocks $R_{j}$, $S_{j}$ which are present within the region delimited by $J_{\ell}$ and $J_{r}$ and
\item the sum $\mathcal R_{13}(t)$ of the strengths of 1- and 3-rarefactions which are present within the region delimited by $J_{\ell}$ and $J_{r}$.
\end{itemize}
Along the lines of~\cite[Lemma~5.2]{CaravennaSpinolo:counterex} one can arrive to the estimates, for a $K>0$,
\[
   \omega^{j+1}\sqrt{\varepsilon} /K  \leq \mathcal |\mathcal V_{13}(t)-\mathcal R_{13}(t) | \leq K \omega^{j+1},
   \qquad
   |\mathcal R_{13}(t)|\leq K \omega^{j}r
 \] 
 at all times when the reflecting waves interacted $j$-th times with $J_{\ell}$ or $J_{r}$.
This is due to the fact that by geometric constraints for being in the triangle with edges $J_{\ell}$ or $J_{r}$ all 1- and 3- rarefactions must be reflected the same number of times as the reflecting shock, or otherwise they merge.
The above estimate imply in particular that the shocks $R_{j}$, $S_{j}$ have strength at least of order $ \omega^{j+1}\sqrt{\varepsilon} /K$.

\vskip.5\baselineskip
This concludes the sketch of stability of the infinite shock structure concerning Theorem~\ref{T:CSa}.

\subsection{Second counterexample}
\label{Ss:secondschaeffer}
Compared to the first construction in \S~\ref{Ss:counterexample}, the main device in order to prevent that rarefactions coming from far destroy the infinite shock structure consists in preventing rarefaction themselves~\cite{CaravennaSpinolo:counterex}. The device is adding to the initial datum a suitable perturbation which is monotone along the eigenvalues: we obtain that all the solutions to the Riemann problems which arise, within a rectangle where we focus the analysis, both at time $t=0$ and at later times, are solved only by shocks, and no rarefaction is present.
The drawback is that the analysis for establishing that only shocks arise at initial time from the more robust initial data is involved, and requires also the new estimates~\cite{CaravennaSpinolo}.

Let us term $\mathfrak S=\mathfrak S(\R;\R^{3})$ the Schwartz space of functions from $\R$ to $\R^{3}$:
\[
\mathfrak S=\left\{V\in C^{\infty}(\R;\R^{3}) \ \Big|\ \sup_{x\in\R} \left(\lvert x\rvert^{\alpha}\lvert D^{\beta} V(x)\rvert\right)<\infty\ \forall \alpha,\beta\in\N\cup\{0\} \right\}.
\]

\begin{theorem}[C., Spinolo-2015]
\label{T:CSb}
Consider the Cauchy problem of~\eqref{E:BaitiJenssen}.
There exists a $C^{\infty}$, compactly supported initial datum $\widetilde U$ with small total variation and a constant $r>0$, depending on $\widetilde U$, such that every initial data belonging to the nonempty ball 
\[
\mathfrak B = \widetilde U+ \left\{ Z\in W^{1,\infty} \text{ t.c.~} \norm{Z}_{W^{1,\infty}}<r \right\}
\]
develops infinitely many shocks in a compact set: it contains a shock pattern like in Figure~\ref{fig:esempio}.
\end{theorem}

Being the topology of the Schwartz space $ \mathfrak S\subset W^{1,\infty}$ 
stronger than the one of $W^{1,\infty}$ itself and since $\widetilde U\in\mathfrak S$, then $\mathfrak B\cap\mathfrak S$ is in turn an open subset of $\mathfrak S$: we eventually furnish a full counterexample to the possibility of extending Schaeffer's or Dafermos' regularity theorems to the class of genuinely nonlinear, strictly hyperbolic systems.



\vskip\baselineskip
\paragraph{\bf Ackowledgments}
The counterexamples of the paper are a joint work with Laura V.~Spinolo (IMATI-CNR, Pavia).
The author is deeply grateful for her precious and friendly collaboration.
Both of us appreciated the lively and constructive atmosphere during HYP2014, which motivated us in improving our result and suggests further developments.
Part of the work was done while the author was affiliated at the OxPDE--University of Oxford, which provided a great working environment.
The author wish to thank also the support of the GNAMPA--INdAM and of the PRIN national project ``Nonlinear Hyperbolic PDE, Dispersive and Transport Equations: theoretical and applicative aspects''.


\end{document}